\documentclass[12pt, twoside]{article}
\usepackage{pcam}

\def\d{{\rm d}}
\def\E{{\rm e}}
\def\i{{\rm i}}
\def\sign{{\rm sign}\,}

\newtheorem{theorem}{Theorem}

\newtheorem{problem}{Problem}
\newtheorem{corollary}{Corollary}

\begin{document}
\English

\title
	[Optimal quadrature formula] 
    {On an optimal quadrature formula for approximation of Fourier integrals in the space $W_2^{(1,0)}$} 

\author
	[Babaev~S.\,S.] 
	{Babaev~S.\,S., Hayotov~A.\,R., Khayriev~U.\,N.} 
    [$^{1,2}$Babaev~S.\,S., $^{2,3}$Hayotov~A.\,R., $^{1}$Khayriev~U.\,N.] 

\email
    {s.s.boboev@buxdu.uz; hayotov@mail.ru; u.n@xayriev@buxdu.uz}

\thanks
 {The work has been done while Samandar  Babaev was visiting Department of Mathematical
Sciences at KAIST, Daejeon, Republic of Korea, fellow of  'El yurt umidi' Foundation of Uzbekistan  for Advanced training course.
}

\organization
    {$^1$ Bukhara State University, 11, M.Ikbol str., Bukhara 200114, Uzbekistan;

     $^2$ V.I.Romanovskiy Institute of Mathematics, Uzbekistan Academy of Sciences,81, M.Ulugbek str., Tashkent 100170, Uzbekistan;

     $^3$ National University of Uzbekistan named after Mirzo Ulugbek, 4, University str., Tashkent 100174, Uzbekistan
     }

\abstract
   {The present paper is devoted to construction of an optimal quadrature formula for approximation of Fourier integrals in the Hilbert space $W_2^{(1,0)}[a,b]$ of non-periodic, complex valued functions.   	

 Here the quadrature sum consists of linear combination of the given function values on uniform grid. The difference between integral and quadrature sum is estimated by the norm of the error functional. The optimal quadrature formula is obtained by minimizing the norm of the error functional with respect to coefficients. In addition, analytic formulas for optimal coefficients are obtained using the discrete analogue of the differential operator $\d^2/\d x^2-1$. Further, the order of convergence of the optimal quadrature formula is studied.
}
 

\titleRus
    [Оптимальная квадратурная формула] 
    {Об оптимальной квадратурной формуле для аппроксимации интегралов Фурье в пространстве $ W_2 ^ {(1,0)} $} 

\authorRus
    [Бабаев ~ С.\,С.] 
    {Бабаев ~ С.\,С., Хаётов ~ А.\,Р., Хайриев ~ У.\,Н.} 
    [$ ^ {1,2} $ Бабаев ~ С.\,С., $ ^ {2,3} $ Хаётов ~ А.\,Р., $ ^ {1} $ Хайриев ~ У.\,Н.] 

\thanksRus
    {Работа была проделана в то время, когда Самандар С. Бабаев посещал Отдел математических наук в KAIST, Тэджон, Республика Корея, в качестве стипендиата Фонда «Эл-юрт умиди».}

\organizationRus
    {$ ^ 1 $ Бухарский государственный университет, ул. М.Икбола, 11, Бухара, 200114, Узбекистан;

     $ ^ 2 $ Институт математики им. В.И.Романовского Академии наук Узбекистана, ул. М.Улугбека, 81, Ташкент 100170, Узбекистан;

$ ^ 3 $ Национальный университет Узбекистана имени Мирзо Улугбека, Ул. Университетская, 4, Ташкент 100174, Узбекистан.
}

\abstractRus
    {Настоящая статья посвящена построению оптимальной квадратурной формулы для аппроксимации интегралов Фурье в гильбертовом пространстве $ W_2 ^ {(1,0)} [a, b] $ непериодических комплекснозначных функций.
  Здесь квадратурная сумма состоит из линейной комбинации значений данной функции на равномерной сетке. Разница между интегралом и квадратурной суммой оценивается по норме функционала погрешности. Оптимальная квадратурная формула получается путем минимизации нормы функционала погрешности по коэффициентам. Кроме того, аналитические формулы для оптимальных коэффициентов получены с использованием дискретного аналога дифференциального оператора $ \ d ^ 2 / \ d x ^ 2-1 $. Также изучена порядок сходимости оптимальной квадратурной формулы.}
\keywordsRus{Оптимальная квадратурная формула, квадратично интегрируемая функция, функционал ошибки, преобразование Фурье, оптимальное приближение}

\UDC{519.644}

\receivedRus{03.02.2019}
\receivedEng{December 3, 2019}

\maketitle

\section{Introduction}
\noindent 
In \cite{BolHayKhu15} and \cite{BolHayMilShad17}, based on Sobolev method, the problem of construction of optimal quadrature formulas
for numerical calculation of Fourier coefficients
\begin{equation}\label{FCoef}
I(\varphi)=\int\limits_0^1\E^{2\pi\i\omega x}\varphi(x)\d x
\end{equation}
with $\omega\in \mathbb{Z}$ was studied in Hilbert spaces $L_2^{(m)}(0,1)$ and $W_2^{(m,m-1)}(0,1)$, respectively. In these works explicit formulas of optimal coefficients were obtained for $m\geq 1$. In particular, for the case $m=1$ the order of convergence of optimal quadrature formulas was studied.

It should be noted that in practice due to the fact that we have discrete values of an integrand the Fourier transforms are reduced
to approximation of integrals of type (\ref{FCoef}) with $\omega\in \mathbb{R}$. For example, the problem of X-ray Computed Tomography (CT) is to reconstruct
the function from its Radon transform. One of the widely used analytic methods of CT reconstruction is the filtered back-projection method
in which the Fourier transforms are used (see \cite{KakSlaney88}).

Recently, in \cite{ZaNov19} authors studied optimal quadrature formulas for arbitrary weighted integrals in the Sobolev space $H^1([0; 1])$. They obtained general formulas for the worst case error depending on the nodes. In particular, for computation of Fourier coefficients of the form (\ref{FCoef})  with $\omega \in \mathbb{R}$ they proved that equidistant nodes are optimal if $n\geq 2.7|\omega|+1$, where $n$ is the number
of the nodes of a quadrature formula.

It should be noted for numerical calculation of integrals (\ref{FCoef}) with real $\omega$ a quadrature formula with explicit coefficients
is required. Therefore, in the present work we study the problem of construction of optimal quadrature formulas in the sense of Sard
for approximate calculation of Fourier integrals of the form (\ref{FCoef}) with $\omega\in \mathbb{R}$ in the space
$W_2^{(1,0)}$.
We obtain explicit formulas for optimal coefficients and calculate the norm of the error functional of the obtained optimal quadrature formula.
We note that the obtained optimal quadrature formula can be used for approximation of Fourier integrals
and reconstruction of a function from its discrete Radon transform.

We note that in  \cite{HaySoom19} it was studied the problem of construction of optimal quadrature
formulas in the sense of Sard for approximate calculation of Fourier integrals of the form (\ref{FCoef}) with $\omega \in \mathbb{R}$ in the space $L_2^{(1)}$.

\section{Construction of optimal quadrature formula}

We examine the following quadrature formula
\begin{equation}\label{eq1}
\int\limits_0^1\E^{2\pi\i\omega x}\varphi(x)\d x\cong \sum\limits_{\beta=0}^NC_{\beta}\varphi(h\beta)
\end{equation}
accompanied by the error
\begin{equation}
\label{eq2}
(\ell,\varphi)=\int\limits_0^1\E^{2\pi \omega \i x}\varphi(x)\d x -\sum\limits_{\beta  = 0}^N {C_\beta  \varphi (h\beta  )}
\end{equation}
where
$$
(\ell,\varphi)=\int\limits_{ - \infty }^\infty  {\ell (x)\varphi(x)\d x}
$$
and respective error functional
\begin{equation}\label{eq3}
\ell (x) =\E^{2\pi \i\omega x}\varepsilon _{[0,1]}(x)-\sum\limits_{\beta  = 0}^N
{C_\beta} \delta (x - h\beta).
\end{equation}
Here $\i^2=-1$, $\omega\in \mathbb{R}$ with $\omega\neq 0$, $C_{\beta}$ are coefficients of the formula (\ref{eq1}), $h=1/N$, $N\in \mathbb{N}$,  $\varepsilon_{[0,1]}(x)$ is the characteristic function of the interval $[0,1]$, and $\delta$
is  Dirac's delta-function.
Function $\varphi$  belongs to the space $W_2^{(1,0)}[0,1]$ which is defined as
\[
\begin{array}{lll}
W_2^{(1,0)}[0,1]&=&\{\varphi:[0,1]\to \mathbb{C}\ |\
\varphi \mbox{ is abs.cont  } \\
& &\mbox{ and }\varphi'\in L_2[0,1]\},
\end{array}
\]
the Hilbert space of complex valued functions.

 The inner product in this space is defined by the equality
\begin{equation}\label{eq4}
\langle\varphi,\psi\rangle =\int\limits_0^1(\varphi'(x)+\varphi(x))(\bar\psi'(x)+\bar\psi(x))\d x,
\end{equation}
where $\bar\psi$ is the complex conjugate function to the function $\psi$ and
the semi-norm of the function $\varphi$ is defined by the
formula
\begin{equation}
 \left\| {\varphi |W_2^{(1,0)}(0,1)} \right\| = \langle\varphi,\varphi\rangle^{\frac{1}{2}} ,\label{eqno(1.3)}
\end{equation}
and   $\int\limits_0^1(\varphi'(x)+\varphi(x))(\bar\varphi'(x)+\bar\varphi(x))\d x<\infty.$

We note that the coefficients $C_{\beta}$ of the formula (\ref{eq1}) depend on $\omega$ and $h$, i.e., $C_{\beta}=C_{\beta}(\omega,h)$.

The error (\ref{eq2}) of the quadrature formula (\ref{eq1}) is a linear functional
in $W_2^{(1,0)*}[0,1]$, where $W_2^{(1,0)*}[0,1]$ is the
conjugate space to the space $W_2^{(1,0)}[0,1]$.

The norm of the error functional (\ref{eq3}) defined as
\begin{equation}\label{eq5}
\left\| \ell \right\|_{W_2^{(1,0)*}[0,1]} = \mathop {\sup
}\limits_{\left\| {\varphi } \right\|_{W_2^{(1,0)}[0,1]} = 1}
\left| {\left( {\ell,\varphi } \right)} \right|
\end{equation}

By the Cauchy-Schwarz inequality the absolute value of the error (\ref{eq2}) is estimated as
$$
|(\ell,\varphi)|\leq \|\varphi\|_{W_2^{(1,0)}[0,1]}\cdot
\|\ell\|_{W_2^{(1,0)*}[0,1]},
$$

The problem of construction of the optimal quadrature formula (\ref{eq1}) in the sense of Sard \cite{Sard}
consists in finding the minimum of the norm (\ref{eq5}) of the error
functional $\ell$ by coefficients $C_{\beta}$ when the nodes
are fixed. We note that here distances between neighbor nodes of the formula (\ref{eq1}) are equal.
This problem, for the quadrature formulas of the form (\ref{eq1}) with $\omega=0$, was first studied in the space $W_2^{(m,m-1)}$ for some $m$ in \cite{ShHay12}, where $W_2^{(m,m-1)}$ is the space of real-valued functions which are square integrable with $m$th generalized derivative.

Therefore, for constructing optimal quadrature formulas of the form (\ref{eq1}) in the sense of Sard
in the space $W_2^{(1,0)}[0,1]$ we need to solve the following problem.

\begin{problem}\label{Prob1}
Find the coefficients $\mathring{C}_\beta $ that give $\inf$ value to $\left\| \ell \right\|_{W_2^{(1,0)*}[0,1]}$, and calculate
\begin{equation}
\left\| {\mathring{\ell}}\right\|_{W_2^{(1,0)*}[0,1]} = \mathop {\inf }\limits_{C_{_\beta  } } \left\| \ell \right\|_{W_2^{(1,0)*}[0,1]}.\nonumber
\end{equation}
\end{problem}

In this section we solve Problem \ref{Prob1} for the cases $\omega\in \mathbb{R}$ with  $\omega\neq 0$ by finding
the norm (\ref{eq5}) and minimizing it by coefficients $C_{\beta}$.

\subsection{The norm of the error functional (\ref{eq3})}

For finding the norm (\ref{eq5}) we use \emph{the extremal function} for the error functional $\ell$ (see, \cite{Sobolev74,SobVas})
which satisfies the following equality
\begin{equation}\label{eq7}
 \left({\ell,\psi _\ell} \right) = \left\| \ell \right\|_{W_2^{(1,0)*}[0,1]} \cdot
 \left\|\psi_\ell \right\|_{W_2^{(1,0)}[0,1]}.
\end{equation}

Since $W_2^{(1,0)}[0,1]$ is a Hilbert space, then the
extremal function $\psi_\ell$ in this space is found by using the Riesz theorem
on the general form of a linear continuous functional on a
Hilbert space. Then for the
functional $\ell$ and for any $\varphi\in
W_2^{(1,0)}[0,1]$ there exists the function $\psi_\ell \in
W_2^{(1,0)}[0,1]$ for which the following equation holds
\begin{equation}
\left( {\ell,\varphi} \right) = \left\langle {\psi _\ell,\varphi}
\right\rangle . \label{eq8}
\end{equation}
Furthermore the equality  $\|\ell\|_{W_2^{(1,0)*}[0,1]} = \|\psi_\ell\|_{W_2^{(1,0)}[0,1]}$ is fulfilled.
Then from (\ref{eq7}), we get
\begin{equation}\label{eq9}
\left(\ell,\psi_{\ell}\right) = \|\ell\|_{W_2^{(1,0)*}[0,1]}^2.
\end{equation}

For  the error functional (\ref{eq3}) to be defined on the space $W_2^{(1,0)}[0,1]$ it should be imposed the following condition
\begin{equation}\label{eq10}
(\ell,\E^{-x})=0,
\end{equation}
which means that the quadrature formula (\ref{eq1}) is exact for $\E^{-x}$.

From (\ref{eq8}) for the extremal function $\psi_{\ell}$ we have the following boundary value problem
\begin{eqnarray}
&&\bar\psi _\ell''(x)-\bar\psi _\ell(x)=- \ell(x),         \label{eq11}\\
&&(\bar\psi _\ell'(x)+ \bar\psi _\ell(x))|^{x=1}_{x=0} =0, \nonumber
\end{eqnarray}
where $\bar\ell$ is the complex conjugate to $\ell$. Then the following holds

\begin{theorem}\label{Thm1}
The solution of the boundary value
problem (\ref{eq11}) is the extremal
function $\psi_{\ell}$ of the error functional $\ell$ and is expressed as
 \begin{equation}\label{eq13}
\psi_\ell(x)=-\bar\ell (x)*G(x)+ \bar d\E^{-x},
 \end{equation}
where
 \begin{equation}\label{eq14}
G(x) = \frac{\sign{x}}{2}\sinh(x),
\end{equation}
$\bar d=d^R+\i d^I$, a complex number, and $*$ is the operation of convolution.
\end{theorem}

Next, we assume that
\begin{equation}\label{eq15}
C_{\beta}=C_{\beta}^R+\i C_{\beta}^I,
\end{equation}
where $C_{\beta}^R$ and $C_{\beta}^I$ are real numbers.
Then, keeping (\ref{eq9}) in mind and using (\ref{eq10}) and (\ref{eq13}) for the norm of the error functional $\ell$ we have
\begin{eqnarray*}
\|\ell\|^2&=&(\ell,\psi_{\ell})\nonumber \\
&=&\int\limits_{-\infty}^{\infty}\ell(x)\psi_{\ell}(x)\d x=-\int\limits_{-\infty}^{\infty}\ell(x)\cdot (\bar \ell(x)*G(x))\d x.
\end{eqnarray*}
Hence, taking (\ref{eq15}) into account, by direct calculation we get
\begin{eqnarray}
\|\ell\|^2&=&-\Bigg[\sum\limits_{\beta=0}^N\sum\limits_{\gamma=0}^N(C_{\beta}^RC_{\gamma}^R+C_{\beta}^IC_{\gamma}^I)\ G(h\beta-h\gamma)\nonumber\\
&&-2\sum\limits_{\beta=0}^NC_{\beta}^R\int\limits_0^1\cos 2\pi \omega x\cdot G(x-h\beta)\d x\nonumber \\
&&-2\sum\limits_{\beta=0}^NC_{\beta}^I\int\limits_0^1\sin 2\pi \omega x\cdot G(x-h\beta)\d x\nonumber\\
&&+\int\limits_0^1\int\limits_0^1\cos[2\pi \omega(x-y)]\cdot G(x-y)\d x\d y\Bigg].\label{eq16}
\end{eqnarray}
Then from (\ref{eq10}), keeping (\ref{eq15}) in mind, we obtain the following equalities
\begin{eqnarray}
&&\sum\limits_{\beta=0}^NC_{\beta}^R\E^{-h\beta}=\int\limits_0^1\cos (2\pi\omega x)\cdot\E^{-x}\ \d x,\label{eq17}\\
&&\sum\limits_{\beta=0}^NC_{\beta}^I\E^{-h\beta}=\int\limits_0^1\sin (2\pi\omega x)\cdot\E^{-x}\ \d x.\label{eq18}
\end{eqnarray}

Thus, we have got the expression (\ref{eq16}) for the norm of the error functional (\ref{eq3}).

Further, in the next section we solve Problem \ref{Prob1}.

\subsection{Minimization of the expression (\ref{eq16}) by coefficients $C_\beta$}

Problem \ref{Prob1} is equivalent to the problem of minimization of the expression
(\ref{eq16}) by $C_\beta^R$ and $C_\beta^I$ using Lagrange method under the conditions (\ref{eq17})-(\ref{eq18}).

Now we consider the function
\begin{eqnarray*}
&&\Psi(C_0^R,C_1^R,...,C_N^R,C_0^I,C_1^I,...,C_N^I,d^R,d^I)\\
&&=\left\| \ell\right\|^2+2d^R\left(\int_0^1\cos (2\pi \omega x)\cdot\E^{-x} \d x-\sum_{\beta=0}^N
C_{\beta}^R\cdot\E^{-h\beta}\right)\\
&&+2d^I\left(\int_0^1\sin (2\pi \omega x)\cdot\E^{-x} \d x-\sum_{\beta=0}^N
C_{\beta}^I\cdot\E^{-h\beta}\right).
\end{eqnarray*}
Equating to zero the partial derivatives of
$\Psi$ by $C_\beta^R$, $C_{\beta}^I$, $(\beta =\overline{0,N})$,
$d^R$ and $d^I$, we get the following
system of linear equations
\begin{eqnarray}
&&\sum\limits_{\gamma=0}^N C_\gamma^R G(h\beta   - h\gamma) +
d^R\E^{-h\beta}= \int\limits_0^1\cos (2\pi \omega x) G(x-h\beta)\d x, \beta=0,...,N,\label{eq19}\\
&&\sum\limits_{\gamma=0}^N C_\gamma^R\E^{-h\gamma}
=\int\limits_0^1\cos 2\pi \omega x\cdot\E^{-x}\ \d x,\label{eq20}\\
&&\sum\limits_{\gamma=0}^N C_\gamma^I G(h\beta   - h\gamma) +
d^I\E^{-h\beta} = \int\limits_0^1\sin 2\pi \omega x G(x-h\beta)\d x, \beta=0,...,N,\label{eq21} \\
&&\sum\limits_{\gamma=0}^N C_\gamma^I\E^{-h\gamma}
=\int\limits_0^1\sin 2\pi \omega x\cdot\E^{-x}\ \d x.\label{eq22}
\end{eqnarray}
Now multiplying both sides of (\ref{eq21}) and (\ref{eq22}) by $\i$ and adding to both sides of
(\ref{eq19}) and (\ref{eq20}), respectively, we get the following system
of $(N+2)$ linear equations with $(N+2)$ unknowns $C_{\gamma}$, $\gamma=0,1,...,N$ and $d$:
\begin{eqnarray}
&&\sum\limits_{\gamma=0}^N C_\gamma G(h\beta   - h\gamma) +
 d\E^{-h\beta} = \int_0^1\E^{2\pi \i\omega x} G(x - h\beta)\d x,\ \beta=0,...,N,\label{eq23} \\
&&\sum\limits_{\gamma=0}^N C_\gamma\cdot\E^{-h\gamma}=\frac{\E^{2\pi \i \omega -1}-1}{2\pi \i \omega -1},\label{eq24}
\end{eqnarray}
where $G(x)$ is defined by (\ref{eq14}).

The system (\ref{eq23})-(\ref{eq24}) has a unique solution. The uniqueness of the solution
of this system can be proved as uniqueness of the solution of the system (3.1)-(3.2) of the work \cite{ShHay12}.
The solution of the system (\ref{eq23})-(\ref{eq24}) gives the minimum to $\left\|\ell\right\|^2 $
in certain value of $C_{\beta}=\mathring{C}_{\beta}$. The quadrature formula of the form (\ref{eq1}) with coefficients
$\mathring{C}_{\beta}$ is called \emph{the optimal quadrature formula} in the sense of Sard and
$\mathring{C}_{\beta}$ are said to be \emph{the optimal coefficients}.

For the convenience, the optimal coefficients $\mathring{C}_{\beta}$ will be denoted as $C_{\beta}$.

The aim of this section is to get the analytic solution of the system (\ref{eq23})-(\ref{eq24}).
For this we  use the concept of discrete argument functions and
operations. The theory of discrete argument functions is
given in \cite{Sobolev74,SobVas}.

Assume that the nodes $x_\beta$ are equally spaced, i.e.,
$x_\beta=h\beta,$ $h$ is a positive small parameter, and functions $\varphi(x)$ and
$\psi(x)$ are complex-valued and defined on the real line $\mathbb{R}$ or on an interval of $\mathbb{R}$.

The function $\varphi (h\beta )$ is {\it a
function of discrete argument } {if it is given on some set
of integer values of} $\beta$.

Further, we need the discrete analogue $D(h\beta)$ of the differential operator $\d^2/\d x^2-1$ which satisfies the following
equality
\begin{equation}\label{eq25}
D(h\beta)*G(h\beta)=\delta_{\d}(h\beta),
\end{equation}
where $\delta_{\d}(h\beta)=\left\{
\begin{array}{rl}
0,&\ |\beta|\neq 0,\\
1,&\ \beta=0
\end{array}
\right. $ ,
$\delta_{\d}(h\beta)$ is the discrete delta-function.

\begin{theorem}\label{Thm2}
The discrete analogue $D(h\beta)$ of the differential operator $\d^2/\d x^2-1$ satisfying (\ref{eq25})
has the form
\begin{equation}\label{eq26}
D(h\beta)=\frac{1}{1-\E^{2h}}\left\{
\begin{array}{rl}
0,&\ |\beta|\geq 2,\\
-2\E^h,&\ |\beta|=1,\\
2(1+\E^{2h}),& \  \beta=0
\end{array}
\right.
\end{equation}
and satisfies the equalities
\begin{eqnarray}\label{eq27}
&&1.\ \  D(h\beta)*\E^{h\beta}=0;\\ \nonumber
&&2.\ \  D(h\beta)*\E^{-h\beta}=0;\\
&&3.\ \  D(h\beta)*G{(h\beta)}=\delta(h\beta). \nonumber
\end{eqnarray}
\end{theorem}

Now we return to our problem.

We consider the coefficients $C_{\beta}$ as a discrete argument function and assume $C_{\beta}=0$ for $\beta<0$ and $\beta>N$. Then the system (\ref{eq23}) and (\ref{eq24}) can be rewriting as follows:

\begin{eqnarray}
&&C_{\beta}*G(h\beta) +
d\E^{-h\beta} = f(h\beta),\ \ \beta  =0,1,...,N, \label{eq28} \\
&&\sum\limits_{\beta=0}^N C_{\beta}\E^{-h\beta}=g_0,\label{eq29}
\end{eqnarray}
where
\begin{eqnarray}
f(h\beta)&=&\frac{\E^{-h\beta}}{4}\cdot\frac{\E^{2\pi\i\omega+1}+1}{2\pi\i \omega+1}-\frac{\E^{h\beta}}{4}\cdot\frac{\E^{2\pi\i\omega-1}+1}{2\pi\i \omega-1}+\frac{\E^{2\pi\i\omega h\beta}}{(2\pi\i \omega+1)(2\pi\i \omega-1)},\nonumber \\ \label{eq30} \\
g_0 &=&\frac{\E^{2\pi\i \omega-1}-1}{2\pi\i \omega-1},\label{eq31}
\end{eqnarray}
and $G(x)$ is defined by (\ref{eq14}).

Now we have the following problem.

\begin{problem}\label{Prob2}
Find $C_{\beta}$ ($\beta=0,1,..,N$) and $d$ which satisfy the system (\ref{eq28})-(\ref{eq29}) for given $f(h\beta)$ and $g_0$.
\end{problem}

Note that Problem \ref{Prob2} is equivalent to Problem \ref{Prob1}.
The main result of this section is the following.

\begin{theorem}\label{Thm3}
Coefficients of the optimal quadrature formulas of the form (\ref{eq1}) in the sense of Sard for
$\omega\in \mathbb{R}$ with $\omega\neq 0$ in the space $W_2^{(1,0)}[0,1]$ have the form
\begin{equation}\label{eq32}
\begin{array}{rcl}
{C}_0&=&\displaystyle \ \frac{1+\E^{2h}+2\pi\i \omega(\E^{2h}-1)-2\E^{(2\pi\i \omega+1) h}}{(\E^{2h}-1)((2\pi \omega)^2+1)},\\
{C}_\beta&=&\displaystyle \ \frac{2(1+\E^{2h}-2\E^{h}\cos 2\pi\omega h)}{(\E^{2h}-1)((2\pi \omega)^2+1)}\cdot\E^{2\pi\i\omega h\beta},\ \beta=1,2,...,N-1,\\
{C}_N&=&\displaystyle \ \frac{\E^{2\pi\i\omega}(1+\E^{2h}-2\pi\i \omega(\E^{2h}-1)-2\E^{(1-2\pi\i \omega) h})}{(\E^{2h}-1)((2\pi \omega)^2+1)}.
\end{array}
\end{equation}
Furthermore, for the square of the norm of the error functional (\ref{eq3}) of the optimal quadrature formula (\ref{eq1}) on the space
$W_2^{(1,0)*}[0,1]$ the following holds
\begin{equation}\label{eq33}
\begin{array}{ll}
\left\|\mathring{\ell}\right\|_{W_2^{(1,0)*}}^2=\frac{1}{(4\pi^2\omega^2+1)^2}\left(4\pi^2\omega^2+1
-\frac{2(1+\E^{2h}-2\E^h\cos (2\pi\omega h))}{h(\E^{2h}-1)}\right).
\end{array}
\end{equation}
\end{theorem}

\emph{Proof.} We consider the following discrete argument function
\begin{equation}\label{eq34}
u(h\beta)=C_\beta*G(h\beta)+d\E^{-h\beta}.
\end{equation}
Then, taking (\ref{eq25}) and (\ref{eq27}) into account, we have
\begin{equation}\label{eq35}
C_{\beta}=D(h\beta)*u(h\beta).
\end{equation}
For calculating the convolution (\ref{eq35}) we need the representation of the function $u(h\beta)$ for all integer values of $\beta$. From (\ref{eq28}) we have
\begin{equation}\label{eq36}
u(h\beta)=f(h\beta)\mbox{ for }\beta=0,1,...,N.
\end{equation}
Now we should find the representation of $u(h\beta)$ for $\beta<0$ and $\beta>N$.

For $\beta\leq 0$ and $\beta\geq N$, using (\ref{eq14}) and (\ref{eq29}), respectively, we get the following
\begin{equation}\label{eq37}
u(h\beta)=
\left\{
\begin{array}{ll}
-\frac{\E^{h\beta}}{4}\ g_0+\left(\frac{1}{4}\sum\limits_{\gamma=0}^NC_{\gamma}\E^{h\gamma}+d\right)\E^{-h\beta},&\ \beta\leq 0,\\
\frac{\E^{h\beta}}{4}\ g_0+\left(-\frac{1}{4}\sum\limits_{\gamma=0}^NC_{\gamma}\E^{h\gamma}+d\right)\E^{-h\beta},& \ \beta\geq N,
\end{array}
\right.
\end{equation}
where $g_0$ is defined by (\ref{eq31}), $\sum\limits_{\gamma=0}^NC_{\gamma}\E^{h\gamma}$ and $d$ are unknowns. We denote
$$a^{-}=\frac{1}{4}\sum\limits_{\gamma=0}^NC_{\gamma}\E^{h\gamma}+d,$$
$$a^{+}=-\frac{1}{4}\sum\limits_{\gamma=0}^NC_{\gamma}\E^{h\gamma}+d.$$
Then from (\ref{eq37}) when $\beta=0$ and $\beta=N$ for these unknowns we obtain the system of two linear equations
\begin{eqnarray*}
&&a^{-}-\frac{1}{4}g_0=f(0),\\
&&a^{+}\E^{-1}+\frac{\E}{4}g_0=f(1).
\end{eqnarray*}
Hence, solving this system, using (\ref{eq30}) and (\ref{eq31}), we get

$a^{-}=\frac{\E^{2\pi\i\omega+1}-1}{4(2\pi\i\omega+1)}$,\      \

$a^{+}=-\frac{\E^{2\pi\i\omega+1}-1}{4(2\pi\i\omega+1)}.$

Hence
\begin{equation}
d=0,\label{eq38}
\end{equation}
\begin{equation}
\frac{1}{4}\sum\limits_{\gamma=0}^NC_{\gamma}\E^{h\gamma}=\frac{\E^{2\pi\i\omega+1}-1}{4(2\pi\i\omega+1)}.
\label{eq39}
\end{equation}
Keeping (\ref{eq38}) and (\ref{eq39}) in mind, and combining (\ref{eq36}) and (\ref{eq37}) we get
\begin{equation*}
u(h\beta)=
\left\{
\begin{array}{ll}
-\frac{\E^{h\beta}}{4}\cdot\frac{\E^{2\pi\i\omega-1}-1}{2\pi\i\omega-1}+\frac{\E^{-h\beta}}{4}\cdot\frac{\E^{2\pi\i\omega+1}-1}{2\pi\i\omega+1},&\ \beta\leq 0,\\
f(h\beta),&\ 0\leq \beta\leq N,\\
\frac{\E^{h\beta}}{4}\cdot\frac{\E^{2\pi\i\omega-1}-1}{2\pi\i\omega-1}-\frac{\E^{-h\beta}}{4}\cdot\frac{\E^{2\pi\i\omega+1}-1}{2\pi\i\omega+1},& \ \beta\geq N.
\end{array}
\right.
\end{equation*}
Now, taking (\ref{eq26}) and (\ref{eq27}) into account, using the last representation of $u(h\beta)$,
and from (\ref{eq35}) by direct calculation for the optimal coefficients ${C}_\beta,$ $\beta=0,1,...,N$,
we obtain analytic formulas (\ref{eq32}).

Now we go to get (\ref{eq33}). We rewrite the equality (\ref{eq16}) in the
following form
\begin{eqnarray}
\|\mathring{\ell}\|^2&=&-\Bigg[\sum\limits_{\beta=0}^NC_\beta^R\left(\sum\limits_{\gamma=0}^NC_\gamma^R
G(h\beta-h\gamma)- \int\limits_0^1\cos 2\pi \omega x\  G(x-h\beta)\d x\right)\nonumber\\
&&+\sum\limits_{\beta=0}^NC_\beta^I\left(\sum\limits_{\gamma=0}^NC_\gamma^I
G(h\beta-h\gamma)- \int\limits_0^1\sin 2\pi \omega x\  G(x-h\beta)\d x\right)\nonumber\\
&& -\sum\limits_{\beta=0}^NC_\beta^R \int\limits_0^1\cos 2\pi \omega x\  G(x-h\beta)\d x
-\sum\limits_{\beta=0}^NC_\beta^I \int\limits_0^1\sin 2\pi \omega x\  G(x-h\beta) \d x\nonumber\\
&&+\int\limits_0^1\int\limits_0^1\cos[2\pi \omega(x-y)]G(x-y)\d x\d y\Bigg].
\label{eq40}
\end{eqnarray}
Since $d=d^R+\i d^I$, taking (\ref{eq38}) into account, we have
$$
d^R=0\mbox{ and } d^I=0.
$$
Therefore, using these last two equalities, from (\ref{eq19}) and (\ref{eq21}) we get the following equalities
$$
\sum\limits_{\gamma=0}^NC_\gamma^R
G(h\beta-h\gamma)- \int\limits_0^1\cos 2\pi \omega x\  G(x-h\beta)\d x=0,\ \beta=0,...,N
$$
and
$$
\sum\limits_{\gamma=0}^NC_\gamma^I
G(h\beta-h\gamma)- \int\limits_0^1\sin 2\pi \omega x\  G(x-h\beta)\d x=0,\ \beta=0,...,N.
$$
Then the expression (\ref{eq40}) for $\|\mathring{\ell}\|^2$ takes the form
\begin{eqnarray*}
\|\mathring{\ell}\|^2&=&\sum\limits_{\beta=0}^NC_\beta^R \int\limits_0^1\cos 2\pi \omega x\  G(x-h\beta)\d x
+\sum\limits_{\beta=0}^NC_\beta^I \int\limits_0^1\sin 2\pi \omega x\  G(x-h\beta) \d x\nonumber\\
&&-\int\limits_0^1\int\limits_0^1\cos[2\pi \omega(x-y)]G(x-y)\d x\d y.
\end{eqnarray*}
Hence calculating the definite integrals, keeping (\ref{eq15}) in mind and using (\ref{eq32}),
after some simplifications we get (\ref{eq33}).
Theorem \ref{Thm3} is proved. \hfill $\Box$

We note that in Theorem \ref{Thm3} the formulas for the optimal coefficients ${C}_\beta$ are decomposed into two parts -- real and imaginary parts. Therefore from the formulas (\ref{eq32}) of Theorem \ref{Thm3} we get the following results.

\begin{corollary}\label{Cor1}
{Coefficients for the optimal quadrature formula of the form
$$
\int\limits_0^1\cos 2\pi \omega x\cdot \varphi(x)\d x\cong \sum\limits_{\beta=0}^NC_{\beta}^R\varphi(h\beta)
$$
in the sense of Sard
in $W_2^{(1,0)}[0,1]$  for $\omega\in \mathbb{R}$ with $\omega\neq 0$ have the form}

\begin{eqnarray*}
{C}_0^R&=&\displaystyle \ \frac{1+\E^{2h}-2\E^{h}\cos 2\pi \omega h }{(\E^{2h}-1)((2\pi \omega)^2+1)},\\
{C}_\beta ^R&=&\displaystyle \ \frac{2(1+\E^{2h}-2\E^{h}\cos 2\pi\omega h)}{(\E^{2h}-1)((2\pi \omega)^2+1)}\cdot\cos(2\pi\omega h\beta),\ \beta=1,2,...,N-1,\\
{C}_N ^R&=&\displaystyle \ \frac{-2\E^{h}\cos(2\pi\omega+2\pi\omega h)+(1+\E^{2h})\cos 2\pi\omega+2\pi\omega\sin 2\pi \omega(\E^{2h}-1) }{(\E^{2h}-1)((2\pi \omega)^2+1)}.
\end{eqnarray*}
\end{corollary}

\begin{corollary}\label{Cor2}
{Coefficients for the optimal quadrature formula of the form
$$
\int\limits_0^1\sin 2\pi \omega x\cdot \varphi(x)\d x\cong \sum\limits_{\beta=0}^NC_{\beta}^I\varphi(h\beta)
$$
in the sense of Sard
in $W_2^{(1,0)}[0,1]$  for $\omega\in \mathbb{R}$ with $\omega\neq 0$ have the form}
\begin{eqnarray*}
{C}_0^I&=&\displaystyle \ \frac{2\pi\omega(\E^{2h}-1)-2\E^{h}\sin 2\pi \omega h }{(\E^{2h}-1)((2\pi \omega)^2+1)},\\
{C}_\beta ^I&=&\displaystyle \ \frac{2(1+\E^{2h}-2\E^{h}\cos 2\pi\omega h)}{(\E^{2h}-1)((2\pi \omega)^2+1)}\cdot\sin(2\pi\omega h\beta),\ \beta=1,2,...,N-1,\\
{C}_N^I&=&\displaystyle \ \frac{-2\E^{h}\sin(2\pi\omega+2\pi\omega h)+(1+\E^{2h})\sin 2\pi\omega-2\pi\omega\cos 2\pi \omega(\E^{2h}-1) }{(\E^{2h}-1)((2\pi \omega)^2+1)}.
\end{eqnarray*}
\end{corollary}

It is easy to see that for $\omega\to 0$ from Theorem \ref{Thm3} we get the optimal of the trapezoidal
quadrature formula in $W_2^{(1,0)}[0,1]$ \cite{ShHay12}.

\begin{corollary}\label{Cor3}
Coefficients of the optimal quadrature formula of the form
\begin{equation}\label{eq41}
\int\limits_0^1 \varphi(x)\d x\cong \sum\limits_{\beta=0}^NC_{\beta}\varphi(h\beta)
\end{equation}
in the space $W_2^{(1,0)}[0,1]$ have the form
\begin{equation*}
\begin{array}{rcl}
{C}_0&=&\displaystyle \frac{\E^{h}-1}{\E^{h}+1},\\
{C}_\beta&=&\displaystyle \frac{2(\E^{h}-1)}{\E^{h}+1},\ \beta=1,2,...,N-1,\\
{C}_N&=&\displaystyle \frac{\E^{h}-1}{\E^{h}+1}
\end{array}
\end{equation*}
and for the square of the norm of the error functional of the optimal quadrature formula (\ref{eq41}) on the space
$W_2^{(1,0)*}[0,1]$ the following holds
\begin{equation*}
\left\|\mathring{\ell}\right\|_{W_2^{(1,0)*}}^2=1-\frac{2(\E^h-1)}{\E^h+1}.
\end{equation*}
\end{corollary}

{\it Remark 1} It should be noted that for fixed $\omega$ from (\ref{eq33}) we get
$$
\|\mathring{\ell}\|^2=\frac{1}{12}h^2-\frac{4\pi^2 \omega^2+3 }{360}h^4+O(h^6),
$$
i.e., the order of convergence of the optimal quadrature formula of the form (\ref{eq1}) is $O(h)$.

{\it Remark 2} In particular, from Theorem \ref{Thm3} in the case $\omega\in \mathbb{Z}$ with $\omega\neq 0$,
we get the results of \cite{BolHayMilShad17}.

{\it Remark 3} The equality (\ref{eq39}) means that the optimal quadrature formula of the form (\ref{eq1}) with coefficients (\ref{eq32})
is exact to $\varphi(x)=\E^{-x}$ because
$$
\int\limits_0^1\E^{2\pi\i \omega x}\E^{-x}\ \d x=\frac{\E^{2\pi \i \omega-1}-1}{2\pi \i \omega-1}.
$$
The equalities (\ref{eq39}) and (\ref{eq29}) provide exactness of our optimal quadrature formula to $\E^{x}$ and $\E^{-x}$ respectively.

\section{Optimal quadrature formulas for the interval [a,b]}

Here from the results of the previous section by a linear transform we get optimal quadrature formulas for the interval $[a,b]$.

We consider construction of optimal quadrature formula of the form
\begin{equation}\label{eq3.1}
\int\limits_a^b\E^{2\pi \i \omega x}\varphi(x)\ \d x\cong \sum\limits_{\beta=0}^NC_{\beta,\omega}[a,b]\varphi(x_\beta)
\end{equation}
in the Hilbert space
$W_2^{(1,0)}[a,b]$. Here $C_{\beta,\omega}[a,b]$ are coefficients and $x_\beta=h\beta+a$ $(\in [a,b])$ are nodes of the formula (\ref{eq3.1}),
$\omega\in \mathbb{R}$, $\i^2=-1$, $h=\frac{b-a}{N}$, $N\in \mathbb{N}$.

Now by linear transformation $x=(b-a)y+a$, where $0\leq y\leq 1$, we get
$$
\int\limits_a^b\E^{2\pi \i \omega x}\varphi(x)\ \d x=(b-a)\E^{2\pi \i \omega a}\int\limits_0^1\E^{2\pi \i \omega (b-a)y}\varphi((b-a)y+a)\d y.
$$
Finally, applying Theorem \ref{Thm3} and Corollary \ref{Cor3} to the integral on the right-hand side of the last equality, we get the following
main result of the present work.

\begin{theorem}\label{Thm4}
Coefficients of the optimal quadrature formula of the form
\begin{equation}\label{eq3.2}
\int\limits_a^b\E^{2\pi \i \omega x}\varphi(x)\ \d x\cong \sum\limits_{\beta=0}^NC_{\beta,\omega}[a,b] \varphi(h\beta+a)
\end{equation}
in the sense of Sard in the space $W_2^{(1,0)}[a,b]$ for
$\omega\in \mathbb{R}$ with $\omega\neq 0$ have the form
\begin{equation}\label{eq3.3}
\begin{array}{rcl}
{C}_{0,\omega}[a,b]&=&\displaystyle (b-a)\E^{2\pi\i\omega a}\ \frac{1+\E^{\frac{2h}{b-a}}+2\pi\i\omega (b-a)(\E^{\frac{2h}{b-a}}-1)-2\E^{\frac{h}{b-a}}\E^{2\pi\i\omega h}}{(4\pi^2 \omega^2(b-a)^2+1)(\E^{\frac{2h}{b-a}}-1)},\\
{C}_{\beta,\omega}[a,b]&=&\displaystyle \ \frac{2(b-a)(1+\E^{\frac{2h}{b-a}}-2\E^{\frac{h}{b-a}}\cos 2\pi\omega h)}{(4\pi^2 \omega^2(b-a)^2+1)(\E^{\frac{2h}{b-a}}-1)}\cdot\E^{2\pi\i \omega (h\beta+a)}\ ,\ \beta=1,2,...,N-1,\\
{C}_{N,\omega}[a,b]&=&\displaystyle (b-a)\E^{2\pi\i\omega b}\ \frac{1+\E^{\frac{2h}{b-a}}-2\pi\i\omega (b-a)(\E^{\frac{2h}{b-a}}-1)-2\E^{\frac{h}{b-a}-2\pi\i\omega h}}{(4\pi^2 \omega^2(b-a)^2+1)(\E^{\frac{2h}{b-a}}-1)},
\end{array}
\end{equation}
and for $\omega=0$ take the form
\begin{equation}\label{eq3.4}
\begin{array}{rcl}
{C}_{0,0}[a,b]&=&\displaystyle (b-a)\frac{\E^{\frac{h}{b-a}}-1}{\E^{\frac{h}{b-a}}+1},\\
{C}_{\beta,0}[a,b]&=&\displaystyle 2(b-a)\frac{\E^{\frac{h}{b-a}}-1}{\E^{\frac{h}{b-a}}+1},\ \beta=1,2,...,N-1,\\
{C}_{N,0}[a,b]&=&\displaystyle (b-a)\frac{\E^{\frac{h}{b-a}}-1}{\E^{\frac{h}{b-a}}+1},
\end{array}
\end{equation}
where $h=\frac{b-a}{N}$.
\end{theorem}

\section{Conclution}

Here for approximation of Fourier integrals in the space $W_2^{(1,0)}[0,1]$
the optimal quadrature formula in the sense of Sard is constructed.
By linear transformation the results are extended  to the case of arbitrary interval $[a,b]$.
That is, for approximation of Fourier integrals in the space $W_2^{(1,0)}[a,b]$ the optimal quadrature formula is obtained.
The obtained optimal quadrature formula can be applied to approximate reconstruction of Computed Tomography images from projections.

\maketitleSecondary

\end{document}